\def\vers{Nov.~25, 2010, v.3}
\magnification=1200
\hsize=6.5truein
\vsize=8.9truein
\font\bigfont=cmr10 at 14pt
\font\mfont=cmr9
\font\sfont=cmr8
\font\mbfont=cmbx9
\font\sifont=cmti8
\def\nin{\noindent}
\def\bs{\bigskip}
\def\ms{\medskip}
\def\sss{\smallskip}
\def\ssb{\raise.2ex\h{${\scriptscriptstyle\bullet}$}}
\def\scirc{\,\raise.2ex\hbox{${\scriptstyle\circ}$}\,}
\def\vsp{\raise10pt\h{$\,$}\raise-6pt\h{$\,$}}
\def\mopl{\hbox{$\bigoplus$}}
\def\mprod{\hbox{$\prod$}}
\def\msum{\hbox{$\sum$}}
\def\mtim{\hbox{$\times$}}
\def\C{{\bf C}}

\def\DD{{\bf D}}
\def\g{\widetilde{\gamma}}
\def\h{\hbox}
\def\H{{\cal H}}
\def\M{{\cal M}}
\def\p{\widetilde{\pi}}
\def\Q{{\bf Q}}
\def\q{\quad}

\def\R{{\bf R}}
\def\S{\widetilde{S}}
\def\SS{{\cal S}}
\def\z{\widetilde{\zeta}}
\def\Z{{\bf Z}}
\def\ZZ{\widetilde{Z}}
\def\Xs{\widetilde{X}_s}
\def\XZ{\widetilde{X}_Z}

\def\supp{{\rm supp}}
\def\End{\hbox{\rm End}}
\def\Ext{\hbox{\rm Ext}}
\def\CH{\hbox{\rm CH}}
\def\Cor{\hbox{\rm Cor}}
\def\Im{\hbox{\rm Im}}
\def\IC{\hbox{\rm IC}}
\def\IH{\hbox{\rm IH}}
\def\Hom{\hbox{\rm Hom}}
\def\Gr{\hbox{\rm Gr}}
\def\MHM{{\rm MHM}}
\def\VHS{{\rm VHS}}
\def\Perv{{\rm Perv}}
\def\Spec{{\rm Spec}}
\def\pR{{}^p\!R}
\def\1{{\hskip1pt}}
\def\simto{\buildrel\sim\over\longrightarrow}
\def\into{\hookrightarrow}
\def\onto{\mathop{\rlap{$\to$}\hskip2pt\hbox{$\to$}}}
\def\({{\rm (}}
\def\){{\rm )}}
\hbox{}
\vskip 1cm

\centerline{\bigfont Relative Chow-K\"unneth decompositions}

\medskip
\centerline{\bigfont for morphisms of threefolds}

\bigskip
\centerline{\it Dedicated to Professor Jacob Murre}

\bigskip
\centerline{Stefan M\"uller-Stach and Morihiko Saito}

\bigskip\medskip
{\narrower\noindent
{\mbfont Abstract.} {\mfont
We show that any nonconstant morphism of a threefold admits a
relative Chow-K\"unneth decomposition.
As a corollary we get sufficient conditions for threefolds
to admit an absolute Chow-K\"unneth decomposition.
In case the image of the morphism is a surface, this implies
another proof of a theorem on the absolute Chow-K\"unneth
decomposition for threefolds satisfying a certain condition,
which was obtained by the first author with P.~L.~del Angel.
In case the image is a curve, this improves in the threefold
case a theorem obtained by the second author where the
singularity of the morphism was assumed isolated and the
condition on the general fiber was stronger.}
\par}

\bigskip\bigskip
\centerline{\bf Introduction}
\footnote{}{{\sifont Date\1}{\sfont:\ \vers}}

\bigskip\noindent
Let $f:X\to S$ be a surjective projective morphism of complex
algebraic varieties with $X$ smooth connected and $S$ reduced.
By the decomposition theorem of Beilinson, Bernstein and Deligne
[BBD], there are noncanonical and canonical isomorphisms
respectively in $D^b_c(S,\Q)$ and $\Perv(S,\Q)$
$$\eqalign{\R f_*\Q_X[\dim X]&
\cong\mopl_i\,\pR^if_*(\Q_X[\dim X])[-i],\cr
\pR^if_*(\Q_X[\dim X])&=\mopl_Z\,\IC_ZL_{Z^o}^i,}
\leqno(0.1)$$
where $Z$ runs over closed irreducible subvarieties of $S$.
Here $\pR^if_*={}^p\H^i\R f_*$ with ${}^p\H_i$ the perverse
cohomology functor, $\Perv(S,\Q)$ is the category of perverse
sheaves on $S$, and $\IC_ZL_{Z^o}^i$ denotes the intersection
complex associated with a $\Q$-local system $L_{Z^o}^i$ defined on
a Zariski-open smooth subvariety $Z^o$ of $Z$, see [BBD].
Moreover, (0.1) holds in the (derived) category of mixed Hodge
modules, and the local system $L_{Z^o}^i$ naturally underlies a
polarizable variation of Hodge structure on $Z^o$ whose weight is
$\dim X-\dim Z+i$, see [Sa1], [Sa2].
Recall that the level of a Hodge structure $H$ is the difference
between the maximal and minimal numbers $p$ such that
$\Gr_F^pH_{\C}\ne 0$, see [D2].
We have the following.

\ms\nin
{\bf Proposition~1.} {\it
Let $n=\dim X$, $m=\dim S$.
Then $L_{Z^o}^i=0$ unless $|i|\le n-m$ with
$Z=S$ or $|i|\le n-\dim Z-2$ with $Z\ne S$.
Moreover, the level of the Hodge structure on each stalk of
$L_{Z^o}^i$ is at most $n-m-|i|$ if $Z=S$, and is at most
$n-\dim Z-|i|-2$ if $Z\ne S$.}

\ms
This is a special case of a more general theorem on the direct
images of arbitrary pure Hodge modules [BS].
Using Proposition~1, we can show the following.

\ms\nin
{\bf Theorem~1.} {\it
Let $f:X\to S$ be a nonconstant projective morphism of
quasi-projective varieties over $\C$.
Assume $X$ is smooth and $3$-dimensional.
Then $f$ admits a relative Chow-K\"unneth decomposition.}

\ms
This means that there are mutually orthogonal projectors
$\pi^f_{i,Z}\in\Cor_S^0(X,X):=\CH_{\dim X}(X\mtim_SX)_{\Q}$ such
that their sum is the diagonal, and their action on the perverse
cohomology sheaf $\pR^jf_*(\Q_X[\dim X])$ is the projection to
$\IC_ZL_{Z^o}^j$ for $j=i$, and vanishes otherwise.
(For the relative Chow-K\"unneth decomposition, see also
[CH1], [GHM1], [NS], [Sa5], etc.)
Note that Theorem~1 in the case $\dim X\le 2$ is rather easy, and
is known to specialists.
For instance, the assertion for $\dim X=2$ follows from
(1.5.1) and Proposition~(1.8) below, see also [Sa5].
In the case $\dim S=1$, a similar assertion on the relative
Chow-K\"unneth decomposition for $f$ was proved in [Sa5]
under the hypothesis that the level of the cohomology of the
generic fibers are at most one and $f$ has at most isolated
singularities.
In this paper we simplify some arguments by using the isomorphisms
in (1.5.1) below, see also [NS], Remark 1.9.

For the moment it is very difficult to prove the Chow-K\"unneth
decomposition for general threefolds $X$.
From Theorem~1 we deduce some sufficient conditions for threefolds
to admit an absolute Chow-K\"unneth decomposition as below.
We will denote by $H^2_{\rm tr}(X,\Q)$ the quotient of $H^2(X,\Q)$
by the subgroup generated by the divisor classes
(called the transcendental part).

\ms\nin
{\bf Theorem~2.} {\it Set $m=\dim S$.
Let $X$ be a smooth complex projective variety of dimension $3$
satisfying the following conditions.

\sss\nin
$(a)$ In case $m=2$, the transcendental part
$H^2_{\rm tr}(X,\Q)$ is generated by the images of
$f^*\IH^2(S,\Q)$ and $H^1(X,\Q)\cup f^*\IH^1(S,\Q)$.

\sss\nin
$(b)$ In case $m=1$, the restriction morphism
$H^2_{\rm tr}(X,\Q)\to H^2_{\rm tr}(X_s,\Q)$ vanishes for general
$s\in S^o$ where $X_s:=f^{-1}(s)$.

\sss\nin
Then $X$ admits an absolute Chow-K\"unneth decomposition.
If moreover $S$ is normal and $f$ has connected fibers
\(replacing $S$ with the Stein factorization if necessary\), then
the absolute projectors $\pi^X_j$ can be obtained by decomposing
the relative projectors $\pi^f_{i,Z}$ in Theorem~$1$.}

\ms
Here $\IH^{\ssb}(S,\Q)$ denotes the intersection cohomology
(see [GM] and also [BBD]).
Since $X$ is smooth, $\IH^{\ssb}(S,\Q)$ in the hypothesis $(a)$
is canonically a direct factor of $H^{\ssb}(X,\Q)$ by the
decomposition theorem~(0.1) together with Proposition~1,
see Remarks~(1.14)(i) below.
Note that $S$ can be replaced by the Stein factorization $S'$ in
Theorem~2, see Remark~(1.14)(iii).
In the case $(a)$ the assertion was shown by the first author with
P.~L.~del Angel ([dAM1], [dAM2]) using another method (by taking a
blow-up of $X$ having a morphism onto a smooth $S$).

Set $r:=\dim X-\dim S=3-m$.
Then $L_{S^o}^0=R^rf_*\Q_X|_{S^o}$, and there is a canonical
decomposition compatible with the variation of Hodge structure
$$L_{S^o}^0=(L_{S^o}^0)^c\oplus(L_{S^o}^0)^{nc},
\leqno(0.2)$$
such that $(L_{S^o}^0)^c$ is constant and
$(L_{S^o}^0)^{nc}$ has no nontrivial global section.
By the global invariant cycle theorem [D2], each stalk of
$(L_{S^o}^0)^c$ coincides with the image of the restriction
morphism by $X_s\into X$ for $s\in S^o$.
Assume $f$ has connected fibers (replacing $S$ if necessary).
The hypothesis $(a)$ then means that
$H^1(S,\IC_S(L_{S^o}^0)^{nc})$ has type $(1,1)$, i.e.
$$F^2H^1(S,\IC_S(L_{S^o,\C}^0)^{nc})=0,
\leqno(0.3)$$
where $(L_{S^o,\C}^0)^{nc}:=(L_{S^o}^0)^{nc}\otimes_{\Q}\C$,
see Remark~(1.14)(iv) below.
Condition $(b)$ means that each stalk of $(L_{S^o}^0)^c$
has type $(1,1)$.
This improves a theorem in [Sa5] for $\dim X=3$.
Note that $H^2(X,\Q)$ is not assumed to be of type $(1,1)$ in our
paper.
Indeed, $H^2(X,\Q)$ contains $H^1(X_s,\Q)^{\rm inv}\otimes
\IH^1(S,\Q)$ in both cases and also $\IH^2(S,\Q)$ in the case
$(a)$, where no conditions are imposed on these Hodge structures.
Here $H^1(X_s,\Q)^{\rm inv}$ denotes the monodromy invariant part.

Part of this work was done during a stay of the second author
at the Mathematical Institute of the University of Mainz.
He would like to thank the staff of the institute for the
hospitality.
Both authors would like to thank J.~Murre for useful discussions,
J.~Nagel for informing us of an error in an earlier version
concerning the proof of (1.5.1), and the referee for useful
comments.
This work is partially supported by Sonderforschungs\-bereich
SFB/TRR 45 (Deutsche Forschungs\-gemeinschaft) and by
Kakenhi 21540037.

\ms\nin
{\bf Convention.} In this paper a variety means a quasi-projective
variety over $\C$.

\bs\bs
\centerline{\bf 1. Preliminaries}

\bs\nin
For the convenience of the reader we give here a short proof of
Proposition~1 using the nearby and vanishing cycle functors [D1].

\ms\nin
{\bf 1.1.~Proof of Proposition~1.}
The assertion can be proved by using the fact that the direct
factors $\IC_ZL_{Z^o}^i$ with $Z\ne S$ are subquotients of
$$\pR^if_*(\varphi_{h,1}\Q_X[\dim X])=
\varphi_{g,1}\pR^if_*(\Q_X[\dim X]),$$
where $h=gf$ with $g$ a function locally defined on $S$ such that
$f$ is smooth over $\{g\ne 0\}$.
Here $\pR^if_*={}^p\H^i\R f_*$ with ${}^p\H_i$ the perverse
cohomology functor (see [BBD]), and $\varphi_{h,1}$ is the unipotent
monodromy part of the vanishing cycle functor [D1] which is shifted
by $-1$ so that it preserves perverse sheaves.
(Note that $\varphi_{g,1}M=M$ if $\supp\,M\subset g^{-1}(0)$.)
We may assume that $h^{-1}(0)$ is a divisor with normal crossings
applying the decomposition theorem to a resolution of singularities.
Then the $\Gr^W_k\varphi_{h,1}\Q_X[\dim X]$ are direct sums of
the constant sheaves supported on intersections of irreducible
components of $h^{-1}(0)$ where $W$ is the monodromy filtration
up to a shift, and moreover the codimensions of the strata are at
least 2 in $X$.
Indeed, $\varphi_{h,1}\Q_X[\dim X]$ is identified with the image
of the logarithm $N$ of the monodromy $T$ on the nearby cycles
$\psi_{h,1}\Q_X[\dim X]$ (see [Sa1], Lemma~5.1.4),
and the assertion for the nearby cycles is rather well-known,
see e.g.\ [Sa2], Th.~3.3 and also [St].
Then Proposition~1 follows by induction on $\dim X$.

\ms\nin
{\bf 1.2.~Relative correspondences.}
Let $X,Y$ be smooth complex quasi-projective varieties with
projective morphisms $f:X\to S$, $g:Y\to S$ where $S$ is a reduced
complex quasi-projective variety.
The group of relative correspondences is defined by
$$\Cor_S^i(X,Y)=\CH_{\dim Y-i}(X\mtim_SY)_{\Q},$$
if $Y$ is equidimensional. In general we take the direct sum over
the connected components of $Y$.
The composition of relative correspondences is defined by using
the pull-back associated to the cartesian diagram (see [Fu])
$$\matrix{X\mtim_SY\mtim_SZ&\to&(X\mtim_SY)\mtim(Y\mtim_SZ)\cr
\downarrow&\vsp&\downarrow\cr
Y&\to&Y\mtim Y,\cr}$$
together with the pushforward by $X\mtim_SY\mtim_SZ\to X\mtim_SZ$,
see [CH1] for details.

We have a canonical morphism
$$\Cor_S^i(X,Y)\to\Cor^i(X,Y):=\Cor_{\C}^i(X,Y),
\leqno(1.2.1)$$
and this is compatible with composition.
So we get a forgetful functor from the category of relative Chow
motives over $S$ to the category of Chow motives over $\C$,
see loc.~cit.

We have moreover the action of correspondences
$$\eqalign{\Cor_S^i(X,Y)
&\to\Hom(\R f_*\Q_X,\R g_*\Q_Y(i)[2i])\cr
&\to\mopl_j\Hom(\pR^jf_*\Q_X,\pR^{j+2i}g_*\Q_Y(i)),\cr}
\leqno(1.2.2)$$
where the decomposition~(0.1) is used for the second morphism.
This action is compatible with the composition of correspondences,
see [CH1], Lemma 2.21.
The isomorphisms in (0.1) can be lifted respectively in $D^b\MHM(S)$
and $\MHM(S)$ where $\MHM(S)$ is the category of mixed Hodge modules,
see [Sa2].
So we can take the first $\Hom$ of (1.2.2) either in $D^b_c(S,\Q)$ or
$D^b\MHM(S)$, and the second $\Hom$ either in $\Perv(S,\Q)$ or
$\MHM(S)$.
Note that the second morphism of (1.2.2) is not an isomorphism
as in (1.5.1) below since the information of morphisms belonging to
higher extension groups $\Ext^i\,(i>0)$ is lost.

\ms\nin
{\bf 1.3.~Relative Chow motives.}
For a projector $\pi\in\Cor_S^0(X,X)$, the relative Chow motive
defined by $\pi$ is denoted by $(X/S,\pi)$.
More precisely, it is an abbreviation of $(X/S,\pi,0)$ in the
usual notation.
We will denote $(X/S,\pi,i)$ by $(X/S,\pi)(i)$ using the notation
of Tate twist.
Morphisms between $(X/S,\pi)(i)$ and $(Y/S,\pi')(j)$ are defined by
$$\Hom((X/S,\pi)(i),(Y/S,\pi')(j))=\pi'\scirc\Cor_S^{j-i}(X,Y)
\scirc\pi.
\leqno(1.3.1)$$
So an isomorphism $(X/S,\pi)\cong(Y/S,\pi')(j)$ as relative
Chow motives is given by morphisms
$$\zeta\in\Cor_S^j(X,Y),\q\zeta'\in\Cor_S^{-j}(Y,X),$$
satisfying the conditions
$$\pi=\zeta'\scirc\zeta,\q\pi'=\zeta\scirc\zeta',
\leqno(1.3.2)$$
together with
$$\zeta=\pi'\scirc\zeta=\zeta\scirc\pi,\q
\zeta'=\pi\scirc\zeta'=\zeta'\scirc\pi'.
\leqno(1.3.3)$$

\ms\nin
{\bf 1.4.~Remark.}
Assume there are $\zeta\in\Cor_S^j(X,Y)$ and
$\zeta'\in\Cor_S^{-j}(Y,X)$ satisfying
$$\pi=\zeta'\scirc\pi'\scirc\zeta,\q
\pi'=\pi'\scirc\zeta\scirc\zeta'\scirc\pi'.
\leqno(1.4.1)$$
Here we assume $\pi'{}^2=\pi'$.
Note that (1.4.1) implies $\pi^2=\pi$.
Replacing $\zeta$ and $\zeta'$ respectively with
$\pi'\scirc\zeta\scirc\pi$ and $\pi\scirc\zeta'\scirc\pi'$,
condition (1.4.1) remains valid using $\pi^2=\pi$ and
$\pi'{}^2=\pi'$.
So (1.3.3) and moreover (1.3.2) are satisfied, since (1.4.1)
becomes (1.3.2) under condition (1.3.3).

\ms\nin
{\bf 1.5.~Case of flat morphisms with relative dimension 1.}
Let $f:X\to S$, $g:Y\to S$ be projective morphisms of complex
quasi-projective varieties such that $X,Y$ are smooth and $f$ is
{\it flat with relative dimension at most} 1.
Then we have isomorphisms (see also [NS], Remark~1.9)
$$\eqalign{\Cor_S^0&(X,Y)\simto\Hom_{D^b\MHM(S)}(\R f_*\Q_X,
\R g_*\Q_Y)\cr
&\simeq\mopl_{i\ge j}\,\Ext_{\MHM(S)}^{i-j}(\pR^if_*\Q_X,
\pR^jg_*\Q_Y).}
\leqno(1.5.1)$$
The last isomorphism follows by lifting the first isomorphism of
(0.1) in $D^b\MHM(S)$ as explained in (1.2).
Applying (1.5.1) to the case $f=g$, the decomposition theorem~(0.1)
implies that $f$ admits a relative Chow-K\"unneth decomposition if
$f$ is flat with relative dimension at most 1.

As for the proof of (1.5.1), let $\DD$ denote the functor
associating the dual, and set $Z=X\mtim_SY$ with canonical
morphisms
$f':Z\to Y$, $g':Z\to X$, $h:Z\to S$.
We have
$$\DD\Q_Y=\Q_Y(\dim Y)[2\dim Y],
\leqno(1.5.2)$$
since $Y$ is smooth.
Consider first the case where the relative dimension is 1.
We have canonical isomorphisms (see [CH1], Lemma 2.21)
$$\eqalign{\Hom&(\R f_*\Q_X,\R g_*\Q_Y)=
\Hom(g^*\R f_*\Q_X,\Q_Y)\cr
=\Hom&(\R f'_*g'{}^*\Q_X,\Q_Y)=\Hom(g'{}^*\Q_X,f'{}^!\Q_Y)\cr
=\Hom&(\Q_Z,(\DD\Q_Z)(1-\dim Z)[2-2\dim Z]).}$$
So (1.5.1) follows from Proposition~(1.15) below.
The argument is similar in case the relative dimension is 0.

\ms\nin
{\bf 1.6.~Orthogonal decomposition of projectors.}
Recall that a projector $\pi'$ is called a {\it direct factor}
(or a refinement) of a projector $\pi$ if
$$\pi'=\pi\scirc\pi'\scirc\pi,\q\h{or equivalently}\q
\pi'=\pi\scirc\pi'=\pi'\scirc\pi.
\leqno(1.6.1)$$
In this case $\pi'':=\pi-\pi'$ is also a projector and is called
the {\it orthogonal complement} of $\pi'$ in $\pi$.
We say that $\pi=\msum_i\,\pi_i$ is an {\it orthogonal decomposition}
of a projector $\pi$ if the $\pi_i$ are mutually orthogonal
projectors.
In this case $\pi_i$ is a direct factor of $\pi$ in the above sense.

\ms\nin
{\bf 1.7.~Good projectors.}
Let $f:X\to S$ be a surjective projective morphism where $X$ is
smooth connected.
We say that a relative projector $\pi\in\Cor_S^0(X,X)$ is a
{\it good projector} if there is a projective morphism $g:Y\to S$
together with a projector $\pi'\in\Cor_S^0(Y,Y)$ and an isomorphism
as relative Chow motives
$$(X/S,\pi)\cong(Y/S,\pi')(-i)\q\h{for some}\,\,\,i\in\Z,$$
such that $\pi'$ is a direct factor of a relative Chow-K\"unneth
projector $\pi^g_{0,Z}$ with $Z:=g(Y)$, and moreover
$$\End(Y/S,\pi')\simto\End_{\MHM(S)}(\M)\subset\End_{\MHM(S)}
(\pR^0g_*(\Q_Y[\dim Y])),
\leqno(1.7.1)$$
where $\M:=\Im\,\pi'\subset\pR^0g_*(\Q_Y[\dim Y])\in\MHM(S)$.

Note that (1.7.1) is satisfied if $Y$ is flat with relative
dimension 1 over $Z$ by (1.5) or if $Y$ is purely 2-dimensional
and is generically finite over $Z$ at each generic point of $Y$
by Proposition~(1.8) below.

Let $\gamma_i\in\Cor_S^0(X,X)$ be mutually orthogonal projectors,
and $\pi\in\Cor_S^0(X,X)$ be a good projector.
Assume $\pi$ is cohomologically orthogonal to the $\gamma_i$ (i.e.\
their actions on the perverse cohomology sheaves are orthogonal).
Set
$$\p=\g\scirc\pi\scirc\g\q\h{with}\,\,\,\,\g:=\mprod_i(1-\gamma_i).
\leqno(1.7.2)$$
Then $\p$ is a projector which is orthogonal to the
$\gamma_i$ using (1.7.1), see also [Sa5].

\ms\nin
{\bf 1.8.~Proposition.} {\it Let $\pi:\S\to S$ be a surjective
projective morphism of purely $2$-dimensional varieties such that
$\S$ is smooth and every irreducible component of $\S$ is dominant
over $S$.
Let $j:U\into S$ be the largest open subset such that the
restriction $\pi_U:\S_U:=\pi^{-1}(U)\to U$ is finite \'etale.
Set $L_U=(\pi_U)_*\Q_{\S_U}$.
Let $s_i$ be the points of $S$ such that $D_i:=\pi^{-1}(s_i)$
has positive dimension.
Let $D_{i,k}$ be the $1$-dimensional irreducible components of $D_i$.
Then $\R\pi_*\Q_{\S}[2]$ is a perverse sheaf naturally underlying a
mixed Hodge module on $S$, and there are canonical isomorphisms
$$\R\pi_*\Q_{\S}[2]=\IC_SL_U\oplus\bigl(\mopl_{i,k}\,
\Q[D_{i,k}]_{s_i}\bigr)\,\,\,\h{in}\,\,\,\Perv(S)\,\,\h{and also in}
\,\,\,\MHM(S),
\leqno(1.8.1)$$
$$\eqalign{&\CH_2(\S\mtim_S\S)_{\Q}\simto\End(\R\pi_*\Q_{\S}[2])\cr
&=\End(\IC_SL_U)\oplus\bigl(\mopl_i\,\End\bigl(\mopl_k\,
\Q[D_{i,k}]_{s_i}\bigr)\bigr).}
\leqno(1.8.2)$$
Here $\Q[D_{i,k}]_{s_i}$ is a sheaf supported
on $\{s_i\}$ and is generated by $[D_{i,k}]_{s_i}$ over $\Q$,
and $\End$ in $(1.8.2)$ can be taken in both $\Perv(S)$ and
$\MHM(S)$.}

\ms\nin
{\it Proof.}
Using the base change theorem for the direct image by a proper
morphism, we can calculate the stalks of the higher direct image
sheaf $(R^2\pi_*\Q_{\S})_s$ for $s\in S$.
It is nonzero if and only if $s=s_i$ for some $i$, and we have
$$(R^2\pi_*\Q_{\S})_{s_i}=\mopl_k\Q[D_{i,k}]_{s_i}.$$
Since $\IC_SL_U=j_{!*}(L_U[2])$, we have the vanishing of
$\H^0\IC_SL_U$.
Then we get (1.8.1) using the decomposition theorem~(0.1).
This implies the last isomorphism of (1.8.2) since there are no
nontrivial morphisms between intersection complexes with different
supports and $L_U$ underlies a variation of Hodge structure of
type $(0,0)$ on $U$.

For the proof of the first isomorphism of (1.8.2), we have
$$\CH_2(\S\mtim_S\S)=\CH^0(\S_U\mtim_U\S_U)\oplus
\bigl(\mopl_i\mopl_{j,k}\,\CH^0(D_{i,j}\mtim D_{i,k})\bigr),$$
since $\dim\S=\dim S=2$.
Then, using the canonical morphism from the associated short exact
sequence
$$0\to\mopl_i\bigl(\mopl_{j,k}\,\CH^0(D_{i,j}\mtim D_{i,k})\bigr)
\to\CH_2(\S\mtim_S\S)\to\CH^0(\S_U\mtim_U\S_U)\to 0,$$
to the corresponding short exact sequence
$$0\to\mopl_i\,\End(\mopl_k\,\Q[D_{i,k}]_{s_i})\to
\End(\R\pi_*\Q_{\S}[2])\to\End(\IC_SL_U)\to 0,$$
the assertion follows.

\ms\nin
{\bf 1.9.~Notation.}
The projector corresponding to $(id,0)$ by Proposition~(1.8) will
be denoted by
$$\pi_{\S_U/U}\in\Cor_S^2(\S,\S)=\CH_2(\S\mtim_S\S)_{\Q}.$$
By definition its action is the identity on $\IC_SL_U$, and
vanishes on the other direct factors.

\ms\nin
{\bf 1.10.~Proposition.} {\it Assume $S$ projective.
With the above notation, let $\pi_j$ denote the absolute
Chow-K\"unneth projectors for $\S$ constructed in {\rm [Mu1]}.
Set $\p_j=\pi_{\S_U/U}\scirc\pi_j\scirc\pi_{\S_U/U}$
as an absolute projector.
Then the $\p_j$ are mutually orthogonal projectors and give
a decomposition of $\pi_{\S_U/U}$ as an absolute projector.}

\ms\nin
{\it Proof.}
It is enough to show that the $\p_j$ for $j\ne 2$ are mutually
orthogonal projectors.
So the assertion is reduced to
$$\pi_i\scirc\pi_{\S_U/U}\scirc\pi_j=\delta_{i,j}\pi_j\q\h{for}
\,\,\,\,i,j\ne 2.
\leqno(1.10.1)$$
Set $\p'=\Delta_{\S}-\pi_{\S_U/U}$ as an absolute projector.
By the direct sum decomposition (1.8.1) its action on the
cohomology $H^j(\S,\Q)$ vanishes for $j\ne 2$, and is the
projection to the subspace generated by the classes of $D_{i,k}$
for $j=2$.
By the construction in [Mu1], the $\pi_j$ for $j\ne 2$ are good
projectors over the base space $S=\Spec\,\C$, see (1.7).
So (1.10.1) for $i=j$ follows since the action of $\pi_{\S_U/U}$
on $H^j(\S,\Q)$ is the identity for $j\ne 2$.
For $i\ne j$, we have to show
$$\pi_i\scirc\p'\scirc\pi_j=0\q\h{if}\,\,\,i,j\ne 2,\,i\ne j.
\leqno(1.10.2)$$
We have isomorphisms of absolute Chow motives
$$(\S,\pi_j)\cong(Y_j,\eta_j)(-k_j)\q(\S,\p')=(Y',\eta')(-1),$$
where $\dim Y_j=0,1,1,0$ and $k_j=0,0,1,2$ for $j=0,1,3,4$
respectively, and $\dim Y'=0$.
So it is enough to show the vanishing of the composition of
morphisms of Chow motives
$$(Y_i,\eta_i)(-k_i)\buildrel{\xi}\over\to(Y',\eta')(-1)
\buildrel{\xi}\over\to(Y_j,\eta_j)(-k_j),
\leqno(1.10.3)$$
assuming that the action of $\xi$ and $\xi'$ on the cohomology
vanishes (since $i,j\ne 2$).
By (1.5.1) for $S=\Spec\,\C$, we can calculate the composition
in the derived category of mixed Hodge structures since $\dim Y'
\le 1$ and $\dim Y'=0$.
Then we get a composition of elements of $\Ext^1$ by the hypothesis
that the action on the cohomology vanishes.
So the composition (1.10.3) vanishes since the higher extension
groups $\Ext^i$ for $i>1$ vanish in the category of mixed Hodge
structures.
This finishes the proof of Proposition~(1.10).

\ms\nin
{\bf 1.11.~Constant part of a relative projector.}
Let $f:X\to S$ be a surjective projective morphism where $X$
is smooth connected and $S$ is reduced.
Let $\pi$ be a relative projector of $X/S$.
Assume $\pi$ has {\it pure relative degree} $j$, i.e.\ the
action of $\pi$ on $\pR^kf_*(\Q_X[\dim X])$ vanishes for $k\ne j$.
We say that $\pi^c$ is the {\it constant part} of $\pi$ if $\pi^c$
is a direct factor of $\pi$, the image of the action of $\pi^c$
on the shifted local system $\pR^jf_*(\Q_X[\dim X])|_{S^o}$ is its
constant part, and there is an isomorphism as relative Chow motives
for some integer $k$:
$$(X/S,\pi^c)\cong(C_S/S,\pi''_S)(-k).$$
Here $C_S:=C\mtim S$ with $C$ an equidimensional smooth projective
variety (which is not assumed to be connected), and
$\pi''_S:=\pi''\mtim[S]\in\Cor_S^0(C_S,C_S)$ is the pull-back of
a direct factor $\pi''$ of the middle Chow-K\"unneth projector
$\pi^C_{\dim C}$ of $C$.
In this case we define $\pi^{nc}:=\pi-\pi^c$, see (1.6).
This is called the {\it nonconstant part} of $\pi$.
We say that the constant part has {\it relative level} $\le i$
if we can take $C$ as above with $\dim C\le i$.

The direct factor $(X/S,\pi^c)$ is well-defined as a relative
Chow motive, if $\pi^c$ has relative level at most 1.
This follows from (1.5.1).
For the well-definedness of $\pi^c$ as a direct factor of $\pi$,
we have to assume, for example, $f$ is flat with relative dimension
$\le 1$ and use (1.5.1).

Let $f:X\to S$ be a surjective projective morphism where $X$
is smooth connected and $S$ is reduced.
We say that a relative projector $\pi$ of $X/S$ has {\it
generically relative level at most $1$}, if there is a dense smooth
open subvariety $U$ of $S$ together with a surjective smooth
projective morphism $g:Y\to U$ and a relative projector $\pi'$ of
$Y/U$ such that $Y$ is equidimensional with relative dimension
$r:=\dim Y-\dim U\le 1$, $\pi'$ is a direct factor of a relative
Chow-K\"unneth projector $\pi^g_{0,S}$ (in particular, the action
of $\pi'$ on $\pR^jg_*(\Q_Y[\dim Y])$ vanishes for $j\ne 0$),
and there is an isomorphism as relative Chow motives over $U$
$$(X/S,\pi)|_U\cong(Y/U,\pi').$$
Here $X|_U:=f^{-1}(U)$ may be assumed smooth over $U$
shrinking $U$ if necessary.
In the above definition we always assume that $\pi^g_{0,S}$ is
induced by (0.1) and (1.5.1) so that we get in the case of relative
dimension $r\le 1$
$$\End((X/S,\pi)|_U)\cong\End(Y/U,\pi')\simto
\End_{\VHS(U)}(M),
\leqno(1.11.1)$$
where $M$ is the image of the action of $\pi'$ on $L:=R^rg_*\Q_Y$
in the category of variations of Hodge structures $\VHS(U)$.

\ms
For the convenience of the reader we show that $\pi^c$ exists if
$\pi$ has generically relative level at most $1$ by using (1.5.1).
For an argument using the theory of abelian schemes (which is
valid also in the positive characteristic case), see [Sa5].

\ms\nin
{\bf 1.12.~Proposition.} {\it
Let $f:X\to S$ be a surjective projective morphism where $X$ is
smooth connected and $S$ is reduced.
In the notation of $(1.11)$, let $\pi$ be a
relative projector of $X/S$ which has generically relative level
at most $1$.
Then the constant part $\pi^c$ of $\pi$ exists.}

\ms\nin
{\it Proof.}
We first consider the case where $g$ in (1.11) has relative
dimension 1. Set
$$L:=R^1g_*\Q_Y\,\,\,\h{in}\,\,\,\VHS(U).$$
Let $L^c$ be the constant part of $L$,
and $L^{nc}$ the orthogonal complement of $L^c$ under a polarization
of $L$ so that
$$L=L^c\oplus L^{nc}\,\,\,\h{in}\,\,\,\VHS(U).$$
By the semisimplicity of $L$, there is no nontrivial morphism
between $L^c$ and $L^{nc}$ in $\VHS(U)$.
Hence the decomposition is compatible with the action of $\pi'$.
In the notation of (1.11.1) we get then
$$M=M^c\oplus M^{nc}\,\,\,\h{in}\,\,\,\VHS(U),$$
where $M^c=M\cap L^c$, $M^{nc}=M\cap L^{nc}$.
This is also compatible with the action of $\pi'$.
So we get by (1.11.1) a canonical orthogonal decomposition of
projectors
$$\pi'=\pi'{}^c+\pi'{}^{nc},$$
corresponding to the above decomposition.
By (1.11.1) there is a smooth projective curve $C$ together
with a projector $\pi''$ which is a direct factor of a
Chow-K\"unneth projector $\pi^C_1$ of $C$ and such that
$$(Y/U,\pi'{}^c)\cong(C_U/U,\pi''_U),$$
where $C_U:=C\mtim U$ and $\pi''_U:=\pi''\mtim[U]$.
Here $C$ can be a general complete intersection in a general fiber
of $g$ using (1.11.1).
Moreover, the above decomposition together with the isomorphism
$(X/S,\pi)|_U\cong(Y/U,\pi')$ implies the orthogonal decomposition
of projectors
$$\pi|_U=(\pi|_U)^c+(\pi|_U)^{nc},$$
such that
$$(X|_U/U,(\pi|_U)^c)\cong(Y/U,\pi'{}^c),\q
(X|_U/U,(\pi|_U)^{nc})\cong(Y/U,\pi'{}^{nc}).$$
So we get
$$(X|_U/U,(\pi|_U)^c)\cong(C_U/U,\pi''_U).$$
By assumption $\pi$ has pure relative degree, say $i$.
Then there are
$$\zeta\in\Cor_U^{-i}(X|_U,C_U),\q
\zeta'\in\Cor_U^i(C_U,X|_U),$$
inducing the above isomorphism, i.e.
$$\zeta'\scirc\zeta=(\pi|_U)^c,\q\zeta\scirc\zeta'=\pi''_U,$$
together with
$$\zeta=(\pi''_U)\scirc\zeta=\zeta\scirc(\pi|_U)^c,\q
\zeta'=(\pi|_U)^c\scirc\zeta'=\zeta'\scirc(\pi''_U).$$
Since $(\pi|_U)^c=(\pi|_U)\scirc(\pi|_U)^c=(\pi|_U)^c\scirc(\pi|_U)$,
the last equalities imply
$$\zeta=\zeta\scirc(\pi|_U),\q\zeta'=(\pi|_U)\scirc\zeta'.$$
Take any extensions
$$\z\in\Cor_S^{-i}(X,C_S),\q\z'\in\Cor_S^i(C_S,X),$$
of $\zeta$ and $\zeta'$ respectively.
Replacing $\z$ and $\z'$ respectively with
$$\pi''_S\scirc\z\scirc\pi\q\h{and}\q
\pi\scirc\z'\scirc\pi''_S,$$
if necessary, we may assume
$$\z=\pi''_S\scirc\z=\z\scirc\pi,\q
\z'=\pi\scirc\z'=\z'\scirc\pi''_S,$$
since the composition of relative correspondences is compatible
with the restriction over $U$.
Using the injection in (1.11.1), we get
$$\pi''_S=\z\scirc\z',$$
since this hold by restricting over $U$. Define
$$\pi^c=\z'\scirc\z.$$
Then $\pi^c$ is a projector and
$$\pi^c=\pi\scirc\pi^c=\pi^c\scirc\pi.$$
So the assertion follows.
The argument is similar in case the relative dimension of $g$ is 0.
This finishes the proof of Proposition~(1.12).

\ms\nin
{\bf 1.13.~Decomposition of a constant relative correspondence.}
With the notation of (1.11), let $\pi^c$ be the constant part of
$\pi$ so that
$$(X/S,\pi^c)\cong(C_S,\pi''_S)(-k).
\leqno(1.13.1)$$
By definition (1.3.1) this isomorphism is induced by
$$\zeta\in\Cor_S^{-k}(X,C_S),\q\zeta'\in\Cor_S^k(C_S,X),$$
satisfying the conditions
$$\zeta'\scirc\zeta=\pi^c,\q\zeta\scirc\zeta'=\pi''_S,$$
together with
$$\zeta=\pi''_S\scirc\zeta=\zeta\scirc\pi^c,\q
\zeta'=\pi^c\scirc\zeta'=\zeta'\scirc\pi''_S.$$
If $S$ admits an absolute Chow-K\"unneth decomposition with
projectors $\pi^S_j$, then we have an orthogonal decomposition
as absolute correspondences
$$\pi^c=\msum_j\,\pi^c_j\q\h{with}\q\pi^c_j=\zeta'\scirc
(\pi''\mtim\pi^S_{j-2k-\dim C})\scirc\zeta,
\leqno(1.13.2)$$
where $k$ is as in (1.13.1), and the action of $\pi^c_j$ on
$H^i(X,\Q)$ vanishes unless $i=j$.

\ms\nin
{\bf 1.14.~Remarks.} (i)
The first isomorphism of (0.1) is not canonical although the
second is.
Set $r=\dim X-\dim S$.
In the notation of (0.1), we have
$$L_{S^o}^{-r}=\H^0f_*\Q_X|_{S^o},$$
by restricting $f$ over $S^o$.
By the decomposition theorem~(0.1), $\IC_SL_{S^o}^{-r}[r]$ is a
direct factor of $\R f_*\Q_X[\dim X]$,
and we have an inclusion morphism
$$\IC_SL_{S^o}^{-r}[r]\into\R f_*\Q_X[\dim X],$$
which splits and induces an injection
$$\IH^j(S,L_{S^o}^{-r})\into H^j(X,\Q).$$
These inclusion morphism are canonical if $\pR^if_*(\Q_X[\dim X])=0$
for $i<-r$, since the negative extension groups vanish,
see also the second isomorphism of (1.5.1).

\sss
(ii) More generally, let $L_k$ denote the filtration on
$\R f_*\Q_X[\dim X]$ defined by the truncation $^p\tau_{\le k}$
(see [BBD]) so that
$$\Gr^L_k(\R f_*\Q_X[\dim X])=\pR^kf_*(\Q_X[\dim X])[-k].$$
Let $L$ denote also the induced filtration on $H^j(X,\Q)$.
There are canonical isomorphisms
$$\Gr^L_kH^j(X,\Q)=H^{j-\dim X-k}(S,\pR^kf_*(\Q_X[\dim X]).$$
This means that the injective morphism
$$H^{j-\dim X-k}(S,\pR^kf_*(\Q_X[\dim X])\into H^j(X,\Q)$$
is canonical modulo $L_{k-1}H^j(X,\Q)$.
This is a generalization of Remark (i) above.

\sss
(iii) Let $f:X\to S$ be as in (0.1).
Let $f':X\to S'$ be the Stein factorization of $f$ with canonical
finite morphism $\pi:S'\to S$ such that $f=\pi\scirc f'$.
We have the decomposition theorem~(0.1) for $f$ and $f'$.
Since the intersection complexes are stable by the direct image
by the finite morphism $\pi$, the direct image by $\pi$ of the
decomposition~(0.1) for $f'$ gives the decomposition~(0.1) for $f$.
This implies that, if the hypothesis (a) or (b) of Theorem~2 is
satisfied, then the same hypothesis holds with $S$ replaced by $S'$,
since $\IC_S\Q$ is a canonical direct factor of $\pi_*\IC_{S'}\Q$
and $f^{-1}(s)$ is the disjoint union of $f'{}^{-1}(s')$ for
$s'\in\pi^{-1}(s)$.

\sss
(iv) Let $f:X\to S$ be as in Theorem~2.
Assume $\dim S=2$ and $f$ has connected fibers.
The decomposition theorem~(0.1) together with Proposition~1
implies
$$H^1(X,\Q)=\IH^1(S,\Q)\oplus H^1(X_s,\Q)^{\rm inv},$$
where $H^1(X_s,\Q)^{\rm inv}$ denotes the monodromy invariant part
which is identified with the stalk of the constant part of
$(L_{S^o}^0)^c$ at $s\in S^o$.
We have
$$\eqalign{\IH^1(S,(L_{S^o}^0)^c)
&=H^1(X_s,\Q)^{\rm inv}\otimes\IH^1(S,\Q)\cr
&=H^1(X,\Q)\cup\IH^1(S,\Q)\,\,\,\h{mod}\,\,\,L_{-1}H^2(X,\Q).}
\leqno(1.14.1)$$
Indeed, the first isomorphism is clear.
For the second isomorphism, note that the cup product is induced by
$$\R f_*\Q_X[3]\otimes\R f_*\Q_X[3]\to(\R f_*\Q_X[3])[3],$$
where the source contains as a direct factor
$$\IC_SL_{S^o}^{-1}[1]\otimes\IC_SL_{S^o}^{-1}[1].$$
This is a perverse sheaf shifted by $4$ since
$\dim\supp\,\H^{-1}(\IC_SL_{S^o}^{-1})=0$ and
$\H^0(\IC_SL_{S^o}^{-1})=0$.
Then its image by the morphism to $(\R f_*\Q_X[3])[3]$ is contained
in $L_{-1}(\R f_*\Q_X[3])[3]$.
So (1.14.1) follows.
Note that we have in the notation of the decomposition~(0.1)
$$\Gr^L_kH^2(X,\Q)=\cases{\IH^2(S,\Q)\oplus
\bigl(\mopl_{s_i}\,L_{\{s_i\}}^{-1}\bigr)&if $\,\,k=-1$,\cr
\IH^1(S,L_{S^0}^0)\oplus\bigl(\mopl_{\dim Z=1}\,
H^0(Z^o,L_{Z^o}^0)\bigr)&if $\,\,k=0$.\cr
\IH^0(S,\Q)(-1)&if $\,\,k=1$.}$$
Here the following Hodge structures have type $(1,1)$
(see Proposition~1)
$$\IH^0(S,\Q)(-1),\q L_{\{s_i\}}^{-1},\q
H^0(Z^o,L_{Z^o}^0)\,\,(\dim Z=1).$$

\sss
(v) In (0.2) and (0.3) in the introduction, it might be possible
to replace $(L_{S^o})^c$ and $(L_{S^o})^{nc}$ respectively with
$(L_{S^o})^{pc}$ and $(L_{S^o})^{npc}$ where $pc$ and $npc$
respectively stand for potentially constant and non-potentially
constant.
The former is defined by the condition that the monodromy group of
the local system is finite (or the local system becomes trivial by
taking the pull-back under a finite \'etale covering
$\rho:S'{}^o\to S^o$).
The latter is the sum of simple local subsystems which are not
potentially constant.
Here we have to use a relative correspondence for $S'/S$ in order
to capture $L^{pc}$ since $\rho_*\rho^*L^{pc}$ is too big.
Moreover, we have to study the relation with the Chow-K\"unneth
projector of $S'$ and the argument is not so simple.
Note that we essentially replace $S$ with the Stein factorization.
In this case it is rather rare that we have $L^c\ne L^{pc}$.
Note also that the hypothesis (a) cannot be stated as in the
form in Theorem~2 if we replace $(L_{S^o,\C})^{nc}$ in (0.3) with
$(L_{S^o},\C)^{npc}$.

\ms
The following is an improvement of [NS], Prop.~1.10, which is needed
for the proof of (1.5.1).
Some argument is similar to 3.5 in later versions (or 2.9 or 2.10
in some earlier versions) of an unpublished preprint [Sa4].

\ms\nin
{\bf 1.15.~Proposition.} {\it Let $X$ be a complex algebraic variety
of dimension at most $d$.
Then we have the bijectivity of the cycle map}
$$\CH_{d-1}(X)_{\Q}\simto\Hom_{D^b\MHM(X)}(\Q_X,(\DD\Q_X)(1-d)[2-2d]).
\leqno(1.15.1)$$

\ms\nin
{\it Proof.}
Set $D={\rm Sing}\,X$ and $U=X\setminus D$.
Let $\pi:X'\to X$ be the normalization.
Set $D'=\pi^{-1}(D)$ with $\pi':D'\to D$ the restriction of $\pi$.
The Chow group $\CH_{d-1}(X)$ does not change by deleting a closed
subvariety of dimension at most $d-2$.
This is the same for the right hand side of (1.15.1)
since $\H^i\DD\Q_X=0$ for $i<-2\dim X$ (see also the proof of [NS],
Prop.~1.10).
Thus we may assume that $X'$, $D$, $D'$ are smooth, $D$ is purely
$(d-1)$-dimensional, and $\pi'$ is \'etale, shrinking $X$ if necessary.

Let $D_i$ be the connected components of $D$,
and $D'_{i,j}$ be the connected components of $\pi^{-1}(D_i)$
with $d_{i,j}$ the degree over $D_i$.
Define
$$E:=\mopl_i\,E_i\q\hbox{with}\q E_i:={\rm Ker}
\bigl(\mopl_j\,d_{i,j}:\mopl_j\Z[D_{i,j}]\to\Z[D_i]\bigr),$$
where $\Z[D_{i,j}]$ is a free $\Z$-module with (formal) generator
$[D_{i,j}]$ (similarly for $\Z[D_i]$), and the morphism
$d_{i,j}:\Z[D_{i,j}]\to\Z[D_i]$ is the multiplication by $d_{i,j}$, which is identified with the trace morphism for $D'_{i,j}\to D_i$.

Let $\SS_i(X)$ be the set of integral (i.e.\ irreducible and reduced)
closed subvarieties of $X$ with dimension $i$.
We have the following commutative diagram of exact sequences
(which is part of the diagram of the snake lemma):
$$\matrix{0&\to&\mopl_{Y'\in\SS_d(X')}\,\C(Y')^*&\simto&
\mopl_{Y\in\SS_d(X)}\,\C(Y)^*&&\cr
\downarrow&&\downarrow&\vsp&\downarrow\cr
E&\into&\mopl_{D'\in\SS_{d-1}(X')}\,\Z[D']&\onto&
\mopl_{D\in\SS_{d-1}(X)}\,\Z[D]\cr
\downarrow&&\downarrow&\vsp&\downarrow\cr
E&\to&\CH_{d-1}(X')&\onto&\CH_{d-1}(X)\cr
\downarrow&&\downarrow&\vsp&\downarrow\cr
0&&0&&0&&}
\leqno(1.15.2)$$
where $\Z[D]$ is a free $\Z$-module with (formal) generator $[D]$,
and the vertical morphism from $\C(Y)^*$ is the divisor map which
associates the multiplicity of a rational function along each
divisor $D$ (and similarly for $Y',D'$).

On the other hand, we have a distinguished triangle in $D^b\MHM(X)$
$$\Q_X\to\pi_*\Q_{X'}\to\mopl_i(\mopl_j\pi'_*\Q_{D'_{i,j}})/\Q_{D_i}
\to.$$
Applying the dual $\DD$ and $\Hom(\Q_X(d-1)[2d-2],*)$, and then
using an isomorphism as in (1.5.2) together with the adjunction
for $\pi^*,\pi_*$, we get an exact sequence
$$E_{\Q}\to\Hom(\Q_{X'},\Q_{X'}(1)[2]))\to
\Hom(\Q_X,\DD\Q_X(1-d)[2-2d]))\to 0,
\leqno(1.15.3)$$
where the $\Hom$ are taken in $D^b\MHM(X')$ or $D^b\MHM(X)$.
Indeed, we have
$$E_{i,{\Q}}=\Hom\bigl(\Q_{D_i},{\rm Ker}\bigl({\rm Tr}:\mopl_j\,
\pi'_*\Q_{D_{i,j}}\to\Q_{D_i}\bigr)\bigr),$$
and the surjectivity of the last morphism of (1.15.3) follows from
$$\Ext^1\bigl(\Q_{D_i},{\rm Ker}\bigl({\rm Tr}:\mopl_j\,
\pi'_*\Q_{D_{i,j}}\to\Q_{D_i}\bigr)\bigr)=0,$$
which is a consequence of the semisimplicity of pure Hodge
modules.

The first morphism of (1.15.3) is induced by the cycle map for
the cycles $[D_{i,j}]$ in $X'$.
Indeed, the latter is induced by the Gysin morphism
$$\Q_{D_{i.j}}\to\Q_{X'}(1)[2],$$
which is the dual of the restriction morphism
$\Q_{X'}\to\Q_{D_{i.j}}$.
Then, comparing (1.15.3) with the third row of (1.15.2) tensored by
$\Q$, the assertion for $X$ is reduced to that for $X'$, and
follows from [21], Prop. 3.4.
This finishes the proof of Proposition~(1.15).

\bs\bs
\centerline{\bf 2. Proof of the main theorems}

\bs\nin
{\bf 2.1.~Proposition.} {\it Let $f:X\to S$ be a nonconstant
surjective projective morphism of complex quasi-projective
varieties where $X$ is smooth, connected, and $3$-dimensional.
Let $Z$ be a closed irreducible subvariety of $S$ and $i$ be
an integer such that $(Z,i)\ne(S,0)$.
Then there is a good projector $\pi^f_{i,Z}\in\Cor_S^0(X,X)$ in the
sense of $(1.7)$ such that its action on the perverse cohomology
sheaf $\pR^jf_*(\Q_X[\dim X])$ is the projection to the direct
factor $\IC_ZL_{Z^o}^j$ for $j=i$, and vanishes for $j\ne i$.}

\ms\nin
{\it Proof.}
By assumption we have $0<\dim S\le \dim X=3$.
We may assume $L_{Z^o}^i\ne 0$ since $\pi^f_{i,Z}=0$ otherwise.
In case $Z=S$, we have $|i|\le\dim X-\dim S$, and
we may assume $\dim S\le 2$ since $(Z,i)\ne(S,0)$.
In case $Z\ne S$, we have $|i|\le 1-\dim Z$ by Proposition~1,
and in particular, $\dim Z\le 1$.

Let $\XZ$ be a desingularization of $X_Z:=f^{-1}(Z)$.
It is denoted by $\Xs$ if $Z=\{s\}$.

\ms\nin
{\bf Case 1} ($\dim Z=0$).
Here $\dim S$ can be arbitrary as in Case 2.
We have $Z=\{s\}$ for some $s\in S$.
Let $\Xs'\subset\Xs$ be the union of the
2-dimensional irreducible components of $\Xs$.
For $|i|\le 1$ we have the injectivity of the composition of
morphisms
$$L_{\{s\}}^i\to \Gr^W_{i+3}H^{i+3}(X_s,\Q)\to H^{i+3}(\Xs,\Q)\to
H^{i+3}(\Xs',\Q).
\leqno(2.1.1)$$
Here the first morphism is induced by the decomposition theorem
(0.1), and is injective.
The second morphism is also injective by the construction of
the weight spectral sequence using a simplicial resolution,
see [D3].
So the composition (2.1.1) is injective for $i=0,1$, since
$\Xs\setminus\Xs'$ is 1-dimensional.
Then the injectivity for $i=-1$ is reduced to the case $i=1$ using
the relative hard Lefschetz theorem for the direct image
$\R f_*\Q_X[\dim X]$ (see [BBD]) since the latter implies an
isomorphism $\eta:L_{\{s\}}^{-1}\simto L_{\{s\}}^1$ where $\eta$ is
the cohomology class of a relative $f$-ample line bundle.

We first consider the case $i=0$.
We take a smooth projective curve $C_s$ over $s$ (which is not
necessarily connected) together with a correspondence
$$\xi\in\Cor^{-1}_S(\Xs',C_s)=\Cor^{-1}(\Xs',C_s)=
\CH^1(\Xs'\mtim C_s)_{\Q},$$
such that the composition below is injective:
$$L_{\{s\}}^0\to H^3(\Xs',\Q)\buildrel{\xi}\over\to H^1(C_s,\Q)(-1),
\leqno(2.1.2)$$
where the first morphism is given by (2.1.1).

Let $\iota'_s:\Xs'\to X$ denote the canonical morphism.
Then (2.1.1) is induced by the morphism of perverse sheaves
$$(\iota'_s)^*:\pR^0f_*(\Q_X[3])\to H^3(\Xs',\Q)_{\{s\}}.$$
This is induced by (1.2.2), and preserves the decomposition by the
support of intersection complexes.
Here $M_{\{s\}}$ for an abelian group $M$ in general denotes
the sheaf supported on $\{s\}$ and whose stalk is $M$.
Thus the injective morphism (2.1.2) is induced by
$$\zeta:=\xi\scirc(\iota'_s)^*\in\Cor_S^{-1}(X,C_s),$$
using (1.2.2) and the compatibility with composition.
Let $\pi^{C_s}_1$ be an absolute Chow-K\"unneth projector for $C_s$.
Replacing $\zeta$ with $\pi^{C_s}_1\scirc\zeta$ if necessary, we may
assume
$$\zeta=\pi^{C_s}_1\scirc\zeta.$$
The dual of (2.1.2) is surjective, and is induced by
${}^t\zeta\in\Cor_S^1(C_s,X)$.
So there is
$$\gamma_{0,s}\in\Cor_S^0(C_s,C_s)=\Cor^0(C_s,C_s),$$
such that
$\gamma_{0,s}=\pi^{C_s}_1\scirc\gamma_{0,s}\scirc\pi^{C_s}_1$
and moreover, setting
$$\pi^f_{0,\{s\}}:=\zeta'\scirc\zeta\in\Cor_S^0(X,X)\q\h{with}\q
\zeta':={}^t\zeta\scirc\gamma_{0,s},$$
the action of $\pi^f_{0,\{s\}}$ on $L_{\{s\}}^0$ is the identity.
(Note that the action on $L_{\{s'\}}^0$ for $s'\ne s$ vanishes by
considering the support.)
The last condition on $\gamma_{0,s}$ depends only on its action
on $H^1(C_s,\Q)$.
By the first condition on $\gamma_{0,s}$ we get
$$\zeta'=\zeta'\scirc\pi^{C_s}_1.$$
Let $H_s\subset H^1(C_s,\Q)$ be the image of (2.1.2) (up to a Tate
twist), and $H_s'$ be its orthogonal complement giving the
orthogonal decomposition
$$H^1(C_s,\Q)=H_s\oplus H_s'.$$
We may assume that the action of $\gamma_{0,s}$ is compatible with
this decomposition and moreover its restriction to $H_s'$ vanishes.

Set
$$\pi'_{\{s\}}:=\zeta\scirc\zeta'=
\zeta\scirc{}^t\zeta\scirc\gamma_{0,s}\in\Cor_S^0(C_s,C_s).$$
The action of $\pi'_{\{s\}}$ on $H^1(C_s,\Q)$ is the projection
to $H_s$ associated to the orthogonal decomposition using the above
hypothesis on the action of $\gamma_{0,s}$.
Since $\pi'_{\{s\}}$ is a direct factor of $\pi^{C_s}_1$ (i.e.
$\pi'_{\{s\}}=\pi^{C_s}_1\scirc\pi'_{\{s\}}\scirc\pi^{C_s}_1$) by
the above argument, it implies that $\pi'_{\{s\}}$ is a projector.
Then $\pi^f_{0,\{s\}}$ is also a projector using the above hypothesis
on $\gamma_{0,s}$ since
$$(\pi^f_{0,\{s\}})^2={}^t\zeta\scirc\pi^{C_s}_1\scirc\gamma_{0,s}
\scirc\pi'_{\{s\}}\scirc\pi^{C_s}_1\scirc\zeta.$$
It is a good projector since (1.7.1) for $\pi'_{\{s\}}$
is clear.
Then, replacing $\zeta$ and $\zeta'$ respectively with
$\pi'_{\{s\}}\scirc\zeta\scirc\pi^f_{0,\{s\}}$ and
$\pi^f_{0,\{s\}}\scirc\zeta'\scirc\pi'_{\{s\}}$ if necessary,
the assertion follows.

The argument is similar for $|i|=1$ where $C$ is replaced by
a disjoint union of a finite number of points.

\ms\nin
{\bf Case 2} ($\dim Z=1,\,Z\ne S$).
Let $\XZ'\subset\XZ$ be the subvariety consisting
the irreducible components whose image in $S$ is $Z$.
Let $\ZZ$ be the Stein factorization of $\XZ'\to Z$.
Note that $\ZZ$ is smooth since $\XZ'$ is smooth and $\dim Z=1$.
Let $\iota'_Z:\XZ'\to X$, $p_Z:\XZ'\to\ZZ$, and
$q_Z:\ZZ\to Z$ denote the canonical morphisms. Set
$$\zeta:=(p_Z)_*\scirc(\iota'_Z)^*.$$
Its action induces an injection
$$L_{Z^o}^0\into (q_Z)_*\Q_{\ZZ}(-1).$$
This is shown by taking a smooth curve intersecting $Z$
transversally at a sufficiently general point of $Z$ and
using the base change over it.
Note that the restriction to the other direct factor of
$\pR^0f_*(\Q_X[\dim X])$ vanishes by the property of the
strict support decomposition.
Then the projector $\pi^f_{0,Z}$ is defined by
$$\pi^f_{0,Z}:=\zeta'\scirc\zeta\in\Cor_S^0(X,X)\q\h{with}\q
\zeta':={}^t\zeta\scirc\gamma_{0,Z},$$
where $\gamma_{0,Z}\in\Cor_S^0(\ZZ,\ZZ)=
\CH^0(\ZZ\mtim_Z\ZZ)_{\Q}$ is chosen so that $\pi^f_{0,Z}$ is
a projector.
As in Case~1 we assume that the action of $\gamma_{0,Z}$ is
compatible with the orthogonal decomposition associated with the
image of the above injective morphism and moreover the action of
$\gamma_{0,Z}$ on the orthogonal complement vanishes.

Set
$$\pi'_Z:=\zeta\scirc\zeta'\in\Cor_S^0(\ZZ,\ZZ).$$
This is also a projector.
Then, replacing $\zeta$ and $\zeta'$ respectively with
$\pi'_Z\scirc\zeta\scirc\pi^f_{0,Z}$ and
$\pi^f_{0,Z}\scirc\zeta'\scirc\pi'_Z$ if necessary,
the assertion follows.

\ms\nin
{\bf Case 3} ($Z=S,\,\dim S=1$).
If $i=-1$, let $Y$ be a sufficiently general hyperplane
section of $X$ which is smooth, 2-dimensional, and flat over $S$.
Let $i_Y:Y\into X$ denote the inclusion.
Then $\zeta$ is defined by $i_Y^*$, and the remaining argument is
similar to the above cases since $Y$ is flat with relative
dimension 1.
The argument is similar for $i=1$ (taking the transpose).

If $i=\pm 2$, we replace $Y$ with a complete intersection
so that $Y$ is smooth, 1-dimensional and flat over $S$.
Then the argument is similar.

\ms\nin
{\bf Case 4} ($Z=S,\,\dim S=2$).
Let $S'\to S$ be the Stein factorization of $f$,
and $\S\to S'$ be a resolution of singularities.
Let $\xi\in\CH_2(X)_{\Q}$ such that $f_*\xi=[S]$.
Taking the pull-back of $\xi$ by $X_{\S}:=X\mtim_S\S\to X$, we get
$$\xi_{\S}\in\Cor_S^0(X,\S)=\CH_2(X_{\S})_{\Q},$$
such that $(f_{\S})_*\xi_{\S}=[\S]$ where $f_{\S}:X_{\S}\to\S$ is
the base change of $f$.
We have
$$[X_{\S}]\in\Cor_S^0(\S,X)=\CH_3(X_{\S})_{\Q}.$$
By Proposition~(1.8) and with the notation of (1.9) we have
$$\pi_{\S_U/U}\scirc\xi_{\S}\scirc[X_{\S}]\scirc\pi_{\S_U/U}=
\pi_{\S_U/U},$$
using the action on the perverse sheaves.
So the projector $\pi^f_{-1,S}$ is defined by
$$\pi^f_{-1,S}=[X_{\S}]\scirc\pi_{\S_U/U}\scirc\xi_{\S},
\leqno(2.1.3)$$
and the assertion follows, see Remark~(1.4).
The argument is similar for $\pi^f_{1,S}$ (taking the transpose).
This finishes the proof of Proposition~(2.1).

\ms\nin
{\bf 2.2.~Proof of Theorem~1.}
The relative projectors $\pi^f_{i,Z}$ for $(Z,i)\ne(S,0)$ are
constructed in Proposition~(2.1).
These can be modified so that they are orthogonal to each other by
(1.7). Then $\pi^f_{0,S}$ is defined to be the remaining so that
$\sum_{i,Z}\pi^f_{i,Z}$ is the diagonal $\Delta_X$,
and the assertion follows.

\ms\nin
{\bf 2.3.~Proof of Theorem~2.}
If the hypothesis (a) or (b) is satisfied, then it is also
satisfied by replacing $S$ with the Stein factorization of $f$,
see Remark~(1.14)(iii).
So we may assume that $S$ is normal and $f$ has connected fibers.
We will decompose every projector $\pi^f_{i,S}$ viewed as an
absolute projector by (1.2.1) into a direct sum of mutually
orthogonal projectors $(\pi^f_{i,S})_j$ in the sense of (1.6)
so that the action of $(\pi^f_{i,S})_j$ on $H^k(X,\Q)$ vanishes
for $k\ne j$.

In the case $Z=\{s_j\}$, the projectors $\pi^f_{i,s_j}$ are
essentially absolute projectors (over $\{s_j\}$), and we do not
have to decompose them further.
So we may assume $\dim Z>0$.
We first consider the case (a) which is more difficult.

\sss\nin
{\bf Case (a.1)}: $\dim Z=1$. We have $i=0$ by Proposition~1, and
Proposition~(1.12) implies the orthogonal decomposition
$$\pi^f_{0,Z}=(\pi^f_{0,Z})^c+(\pi^f_{0,Z})^{nc},$$
since $(X/S,\pi^f_{0,Z})\cong(\ZZ/S,\pi'_Z)$ in the notation of
Case 2 in (2.1).
By (1.13) we have the decomposition
$$(\pi^f_{0,Z})^c=\msum_{j=0}^2(\pi^f_{0,Z})^c_{j+2}.$$
Note that $(\pi^f_{0,Z})^{nc}=(\pi^f_{0,Z})^{nc}_3$, i.e.\
$(\pi^f_{0,Z})^{nc}_j=0$ for $j\ne 3$.

\ms\nin
{\bf Case (a.2)}: $Z=S,\,|i|=1$.
Here Proposition~(1.12) and (1.13) are not sufficient.
Let $\S\to S$ be a resolution of singularities, and
$\pi^{\S}_j\in\Cor^0(\S,\S)_{\Q}$ be the Chow-K\"unneth
projectors for $\S$, see [Mu1].
Here we may assume that $\pi^{\S}_j$ for $j\ne 2$ are good
projectors in the sense of (1.7) over $S=pt$ since only curves are
used in the construction.
By (2.1.3) we have
$$\pi^f_{-1,S}=\zeta'_{-1,S}\scirc\zeta_{-1,S}\q\h{with}$$
$$\zeta'_{-1,S}=\gamma'_{-1,S}\scirc[X_{\S}]\scirc\pi_{\S_U/U},\q
\zeta_{-1,S}=\pi_{\S_U/U}\scirc\xi_{\S}\scirc\gamma_{-1,S}.$$
Here $\gamma_{-1,S}$ and $\gamma'_{-1,S}$ are added in order to
get mutually orthogonal projectors as in (1.7.2), see the proof of
Theorem~1 in (2.2).
By Proposition~(1.10) we have absolute projectors
$$(\pi^f_{-1,S})_j:=\zeta'_{-1,S}\scirc\pi^{\S}_j\scirc
\zeta_{-1,S}\in\Cor^0(X,X),$$
giving the orthogonal decomposition
$$\pi^f_{-1,S}=\msum_{j=0}^4\,(\pi^f_{-1,S})_j\q\h{in}\,\,\,
\Cor^0(X,X).$$
Note that $f^*\IH^{2}(S,\Q)\subset H^{2}(X,\Q)$ coincides with
the image of the action of $(\pi^f_{-1,S})_2$ on $H^2(X,\Q)$ by
the decomposition (1.8.1).
Similarly we get $(\pi^f_{1,S})_{j+2}$ giving the orthogonal
decomposition
$$\pi^f_{1,S}=\msum_{j=0}^4\,(\pi^f_{1,S})_{j+2}\q\h{in}\,\,\,
\Cor^0(X,X).$$

\ms\nin
{\bf Case (a.3)}: $Z=S,\,i=0$.
Proposition~(1.12) implies the orthogonal decomposition
$$\pi^f_{0,S}=(\pi^f_{0,S})^c+(\pi^f_{0,S})^{nc},$$
and we have by (1.13) the decomposition
$$(\pi^f_{0,S})^c=\msum_{j=0}^4(\pi^f_{0,S})^c_{j+1}.$$
We need the assumption (a) to construct the decomposition
$$(\pi^f_{0,S})^{nc}=\msum_{j=1}^3(\pi^f_{0,S})^{nc}_{j+1},$$
Here $(\pi^f_{0,S})^{nc}_{j+1}=0$ for $j=0,4$, since
$H^j(S,\IC_S(L_{S^o}^0)^{nc})=0$ for $|j|=2$
(using the vanishing of $\Gamma(S^o,(L_{S^o}^0)^{nc})$).
By condition (a) the image of the action of $(\pi^f_{0,S})^{nc}$
has type (1,1), see Remark~(1.14)(iii).
So there is an absolute projector $\gamma$ such that the action of
$\gamma$ on $H^2(X,\Q)$ coincides with that of $(\pi^f_{0,S})^{nc}$
and $(X,\gamma)$ is isomorphic to a finite direct sum of copies of
$(pt,id)(-1)$.
Let $\zeta$ and $\zeta'$ be algebraic cycles on a disjoint union
of copies of $X$ inducing the last isomorphism so that
$$\gamma=\zeta'\scirc\zeta,\q\zeta\scirc\zeta'=id.$$
Then we can set
$$(\pi^f_{0,S})^{nc}_2=(\pi^f_{0,S})^{nc}_2\scirc\gamma
\scirc(\pi^f_{0,S})^{nc}_2,$$
since $\zeta\scirc(\pi^f_{0,S})^{nc}\scirc\zeta'=id$ and this
implies that $(\pi^f_{0,S})^{nc}_2$ is an idempotent.

The argument is similar for $(\pi^f_{0,S})^{nc}_4$ where
$H^2(X,\Q)$ and $(pt,id)(-1)$ are respectively replaced by
$H^4(X,\Q)$ and $(pt,id)(-2)$.
The orthogonality of $(\pi^f_{0,S})^{nc}_2$ and
$(\pi^f_{0,S})^{nc}_4$ follows from the fact that the projectors
factors through the direct sums of copies of $(pt,id)(-1)$ or
$(pt,id)(-2)$.
Then $(\pi^f_{0,S})^{nc}_3$ is defined by the remaining so that
the above decomposition holds.

Thus we get an orthogonal decomposition of every $\pi^f_{i,Z}$ as
explained at the beginning of this subsection, and Theorem~2 is
proved in the case (a).

\ms\nin
{\bf Case (b)}. We may assume $Z=S$ by the argument at the
beginning of this subsection.
By (1.12) we have orthogonal decompositions as relative projectors
$$\pi^f_{i,S}=(\pi^f_{i,S})^c+(\pi^f_{i,S})^{nc},$$
such that the image of the action of $(\pi^f_{i,S})^c$ on the local
systems is the constant part of $L^i_{S^o}$, see (1.6).
Here we use condition $(b)$ in the case $i=0$,
since it implies that the invariant part of $H^2(X_s,\Q)$ has
type (1,1) and is generated by divisor classes.
By (1.13) we then get an orthogonal decomposition as absolute
projectors
$$(\pi^f_{i,S})^c=\msum_{j=0}^2\,(\pi^f_{i,S})^c_{i+j+2}.$$
On the other hand, we do not have to decompose further
$(\pi^f_{i,S})^{nc}$, i.e.
$$(\pi^f_{i,S})^{nc}=(\pi^f_{i,S})^{nc}_{i+3},\q
(\pi^f_{i,S})^{nc}_{i+j+2}=0\,\,\,(j\ne 1),$$
since $H^0(S^o,(L_{S^o}^i)^{nc})=0$ and $\dim S=1$.

Thus we get an orthogonal decomposition of every $\pi^f_{i,Z}$ as
explained at the beginning of this subsection, and Theorem~2 is
proved also in the case (b).

\bs\bs
\centerline{{\bf References}}

\ms
{\mfont
\item{[dAM1]}
P.~L.~del Angel and S.~M\"uller-Stach, Motives of uniruled $3$-folds,
Compos.\ Math.\ 112 (1998), 1--16.

\item{[dAM2]}
P.~L.~del Angel and S.~M\"uller-Stach, On Chow motives of $3$-folds,
Transactions of the AMS 352 (2000), 1623--1633.

\item{[BBD]}
A.~Beilinson, J.~Bernstein and P.~Deligne, Faisceaux pervers,
Ast\'erisque, vol. 100, Soc. Math. France, Paris, 1982.

\item{[BS]}
P.~Brosnan and M.~Saito, Decomposition theorem and level of Hodge
modules, preprint.

\item{[CH1]}
A.~Corti and M.~Hanamura, Motivic decomposition and intersection
Chow groups, I, Duke Math. J. 103 (2000), 459--522.

\item{[CH2]}
A.~Corti and M.~Hanamura, Motivic decomposition and intersection
Chow groups II, preprint.

\item{[D1]}
P.~Deligne, Le formalisme des cycles \'evanescents, in SGA7 XIII
and XIV, Lect.\ Notes in Math.\ 340, Springer, Berlin, 1973,
pp.\ 82--115 and 116--164.

\item{[D2]}
P.~Deligne, Th\'eorie de Hodge II, Publ. Math. IHES 40 (1971),
5--58.

\item{[D3]}
P.~Deligne, Th\'eorie de Hodge III, Publ. Math. IHES 44 (1974),
5--77.

\item{[Fu]}
W.~Fulton, Intersection theory, Springer, Berlin, 1984.

\item{[GHM1]}
B.~B.~Gordon, M.~Hanamura and J.~P.~Murre, Relative Chow-K\"unneth
projectors for modular varieties, J. Reine Angew. Math. 558 (2003),
1--14.

\item{[GHM2]}
B.~B.~Gordon, M.~Hanamura and J.~P.~Murre, Absolute Chow-K\"unneth
projectors for modular varieties, J. Reine Angew. Math. 580 (2005),
139--155.

\item{[GM]}
M.~Goresky and R.~MacPherson, Intersection homology theory,
Topology 19 (1980), 135--162.

\item{[Mu1]}
J.~P.~Murre, On the motive of an algebraic surface, J. Reine Angew.
Math. 409 (1990), 190--204.

\item{[Mu2]}
J.~P.~Murre, On a conjectural filtration on Chow groups of an
algebraic variety, Indag. Math. 4 (1993), 177--201.

\item{[NS]}
J.~Nagel and M.~Saito, Relative Chow-K\"unneth decompositions for
conic bundles and Prym varieties, Int.\ Math.\ Res.\ Not.\
2009, no. 16, 2978--3001.

\item{[Sa1]}
M.~Saito, Modules de Hodge polarisables, Publ. RIMS, Kyoto Univ.
24 (1988), 849--995.

\item{[Sa2]}
M.~Saito, Mixed Hodge modules, Publ. RIMS, Kyoto Univ. 26
(1990), 221--333.

\item{[Sa3]}
M.~Saito, Hodge conjecture and mixed motives, I,
Proc.\ Sympos.\ Pure Math., 53, Amer.\ Math.\ Soc., Providence,
RI, 1991, pp. 283--303.

\item{[Sa4]}
M.~Saito, Bloch's conjecture, Deligne cohomology and higher Chow
groups, unpublished preprint (arXiv:9910113v14, 2002).

\item{[Sa5]}
M.~Saito, Chow-K\"unneth decomposition for varieties with low
cohomological level, preprint (math.AG/0604254).

\item{[St]}
J.H.M.~Steenbrink, Limits of Hodge structures, Inv. Math. 31
(1975/76), no. 3, 229--257.

\medskip
{\sfont
\baselineskip=10pt
Mathematisches Institut der Johannes Gutenberg Universit\"at Mainz,
Staudingerweg 9, 55099 Mainz

\sss
RIMS Kyoto University, Kyoto 606-8502 Japan

\smallskip
\vers
}}
\bye